\documentclass[12pt]{amsart}
\usepackage{amssymb,latexsym}
\usepackage{pdfsync}

\newdimen\AAdi%
\newbox\AAbo%
\def\AAk#1#2{\s_etbox\AAbo=\hbox{#2}\AAdi=\wd\AAbo\kern#1\AAdi{}}%
\def\AAr#1#2#3{\s_etbox\AAbo=\hbox{#2}\AAdi=\ht\AAbo\raise#1\AAdi\hbox{#3}}%
\font\tenmsb=msbm10 at 12pt
\font\sevenmsb=msbm7 at 8pt
\font\fivemsb=msbm5 at 6pt
\newfam\msbfam
\textfont\msbfam=\tenmsb
\scriptfont\msbfam=\sevenmsb
\scriptscriptfont\msbfam=\fivemsb
\def\Bbb#1{{\tenmsb\fam\msbfam#1}}

\newcommand{\beq}{\begin{equation}}
\newcommand{\eeq}{\end{equation}}
\newcommand{\beqr}{\begin{eqnarray}}
\newcommand{\eeqr}{\end{eqnarray}}
\newcommand{\ba}{\begin{array}}
\newcommand{\ea}{\end{array}}

\begin{document}

\newtheorem{thm}{Theorem}
\newtheorem{lem}{Lemma}
\newtheorem{cor}{Corollary}
\newtheorem{rem}{Remark}
\newtheorem{pro}{Proposition}
\newtheorem{defi}{Definition}
\newtheorem{conj}[thm]{Conjecture}
\newcommand{\noi}{\noindent}
\newcommand{\dis}{\displaystyle}
\newcommand{\mint}{-\!\!\!\!\!\!\int}
\numberwithin{equation}{section}

\def \bx{\hspace{2.5mm}\rule{2.5mm}{2.5mm}}
\def \vs{\vspace*{0.2cm}}
\def\hs{\hspace*{0.6cm}}
\def \ds{\displaystyle}
\def \p{\partial}
\def \O{\Omega}
\def \o{\omega}
\def \b{\beta}
\def \m{\mu}
\def \l{\lambda}
\def\L{\Lambda}
\def \ul{u_\lambda}
\def \D{\Delta}
\def \d{\delta}
\def \k{\kappa}
\def \s{\sigma}
\def \e{\varepsilon}
\def \a{\alpha}
\def \tf{\tilde{f}}
\def\cqfd{%
\mbox{ }%
\nolinebreak%
\hfill%
\rule{2mm} {2mm}%
\medbreak%
\par%
}
\def \pr {\noindent {\it Proof.} }
\def \rmk {\noindent {\it Remark} }
\def \esp {\hspace{4mm}}
\def \dsp {\hspace{2mm}}
\def \ssp {\hspace{1mm}}

\def\la{\langle}\def\ra{\rangle}

\def \u{u_+^{p^*}}
\def \ui{(u_+)^{p^*+1}}
\def \ul{(u^k)_+^{p^*}}
\def \energy{\int_{\R^n}\u }
\def \sk{\s_k}
\def \mo{\mu_k}
\def\cal{\mathcal}
\def \I{{\cal I}}
\def \J{{\cal J}}
\def \K{{\cal K}}
\def \OM{\overline{M}}

\def\n{\nabla}

\def\fk{{{\cal F}}_k}
\def\M1{{{\cal M}}_1}
\def\Fk{{\cal F}_k}
\def\Fl{{\cal F}_l}
\def\FF{\cal F}
\def\Gk{{\Gamma_k^+}}
\def\n{\nabla}
\def\uuu{{\n ^2 u+du\otimes du-\frac {|\n u|^2} 2 g_0+S_{g_0}}}
\def\uuug{{\n ^2 u+du\otimes du-\frac {|\n u|^2} 2 g+S_{g}}}
\def\sku{\sk\left(\uuu\right)}
\def\qed{\cqfd}
\def\vvv{{\frac{\n ^2 v} v -\frac {|\n v|^2} {2v^2} g_0+S_{g_0}}}
\def\vvs{{\frac{\n ^2 \tilde v} {\tilde v}
 -\frac {|\n \tilde v|^2} {2\tilde v^2} g_{S^n}+S_{g_{S^n}}}}
\def\skv{\sk\left(\vvv\right)}
\def\tr{\hbox{tr}}
\def\pO{\partial \Omega}
\def\dist{\hbox{dist}}
\def\RR{\Bbb R}\def\R{\Bbb R}
\def\C{\Bbb C}
\def\B{\Bbb B}
\def\N{\Bbb N}
\def\Q{\Bbb Q}
\def\Z{\Bbb Z}
\def\PP{\Bbb P}
\def\EE{\Bbb E}
\def\F{\Bbb F}
\def\G{\Bbb G}
\def\H{\Bbb H}
\def\SS{\Bbb S}\def\S{\Bbb S}

\def\div{\hbox{div}\,}

\def\lcf{{locally conformally flat} }

\def\circledwedge{\setbox0=\hbox{$\bigcirc$}\relax \mathbin {\hbox
to0pt{\raise.5pt\hbox to\wd0{\hfil $\wedge$\hfil}\hss}\box0 }}
%\ignorespaces%

\def\sss{\frac{\s_2}{\s_1}}

\date{\today}
%\date{2019-05-03}
\title[ Extensions of Bonnet--Myers' type theorems ]{ Extensions of Bonnet--Myers' type theorems with the Bakry--Emery Ricci curvature }

\author{}

\author[Qiu]{Hongbing Qiu}
\address{School of Mathematics and Statistics\\ Wuhan University\\Wuhan 430072,
China,
and Hubei Key Laboratory of Computational Science \\ Wuhan University\\Wuhan 430072,
China
%, and   Max Planck Institute for Mathematics in the
 % Sciences\\Inselstr. 22\\D-04103 Leipzig, Germany
 }
 \email{hbqiu@whu.edu.cn}

% \author[Wan]{Jianming Wan}
%\address{School of Mathematics\\ Northwest University\\Xi'an 710127,
%China
 %}
 %\email{wanj${\underline{\quad}}$m@aliyun.com}

\begin{abstract}

In this paper, we prove the extensions of Bonnet--Myers' type theorems obtained by Calabi and Cheeger--Gromov--Taylor via Bakry--Emery Ricci curvature, which generalize the results of \cite{FG, Lim1, Wan, Wang, WW, Wu}.

\vskip12pt



\end{abstract}
\maketitle
\section{Introduction}

Let $(M, g)$ be an $n$-dimensional complete Riemannian manifold. The celebrated Bonnet--Myers theorem states that if the Ricci curvature of $M$ has a positive lower bound, then $M$ must be compact. Calabi \cite{Cal} generalized this result as follows\\

\noindent{\bf Theorem}(Calabi) {\it  Let $(M, g)$ be an $n$-dimensional complete Riemannian manifold with nonnegative Ricci curvature. If for some point $p_0 \in M$, every geodesic ray issuing from $p_0$ has the property that 
\[
\limsup_{a\to \infty} \left( \int_0^{a} \left( {\rm Ric}_M(s) \right)^{\frac{1}{2}}ds - \frac{1}{2}\log a  \right) =\infty,
\]
then $M$ is compact. 
 } \\

In particular,  this indicates that if the Ricci curvature of $M$ satisfies 
\[
{\rm Ric}_M(x) \geq \frac{1}{(4-\e)(1+r(x))^2},
\]
where $r(x)=d(p_0, x)$ is the distance function and $\e \in (0, 4)$ is a constant, then $M$ is compact(c.f.\cite{WSY}). Cheeger--Gromov--Taylor \cite{CGT} and Garding \cite{Gar} also proved a similar result. The idea of their proof relies on studying carefully the index form or the second variation.  Recently, Wan \cite{Wan} gave a complementary extension of Calabi and Cheeger--Gromov--Taylor's theorems by showing that the manifold has no ray issuing from some point.

Let $V$ be a smooth vector field on $M$, when Bakry-Emery \cite{BE} studied the diffusion processes, they defined a generalization of the Ricci tensor  of $M$ by 
\[
{\rm Ric}^m_{V} := {\rm Ric} -\frac{1}{2} L_V g - \frac{1}{m-n} V^*\oplus V^*, 
\]
where ${\rm Ric}$ is the Ricci curvature of $M$, $L$ denotes the Lie derivative, $V^*$ is the dual 1-form to $V$ and $m > n$ is a constant. Usually, ${\rm Ric}^m_{V}$ is called the $m$-Bakry--Emery Ricci tensor. When $m=\infty$, one denotes 
\[
{\rm Ric}_V := {\rm Ric} -\frac{1}{2} L_V g.
\]
This is also referred to as the $\infty$-Bakry-Emery Ricci tensor.  In particular, if 
$
{\rm Ric}_V:= \rho g
$
for some constant $\rho$, then $(M, g)$ is called a Ricci soliton, which was introduced by Hamilton \cite{Ham1}. It is well known that Ricci solitons often arise as limits of dilations of singularities in the Ricci flow \cite{Cao1, CZ, Ham2,  Sesum}. When $V = -\n f$ for some smooth function $f$ on $M$, then the Bakry-Emery Ricci tensor ${\rm Ric}_V$ becomes  
\[
{\rm Ric}_f := {\rm Ric} + {\rm Hess}(f).
\]
This tensor is closely related with the smooth metric measure space $(M, g, e^{-f}d{\rm vol}_g)$, which is especially relevant for collapsing (see e.g. \cite{CC, Fuk, Gro}) and plays an important role in Perelman's work on Ricci flow.

In the past two decades, the previous  Bonnet--Myers' type theorems were generalized to the case that the complete Riemannian manifold is equipped with Bakry--Emery Ricci cuvature under different  conditions (see e.g. \cite{BQ, COS}, \cite{HL}--\cite{PS}, \cite{Soy}--\cite{Tad3},% LH, Li, Lim1, Lim2, Lot, MRV, Mor, PS, Soy, Tad1, Tad2, Tad3,
 \cite{Wang, WW, Wu, Zhang}).   
 %As in \cite{CJQ, WW}, we define the Bakry--Emery Ricci curvature

 %Throughout this paper, we use the convention that $m=n$ if and only if $V\equiv 0$. 

Fern\'andez-L\'opez and Garc\'ia-R\'io \cite{FG} proved that if the Bakry--Emery Ricci curvature ${\rm Ric}_V$ of a complete Riemannian manifold $(M, g)$ has a positive lower bound, and $|V|$ is bounded, then $M$ is compact. Later, by applying the mean curvature comparison theorem to the excess function, Wei--Wylie \cite{WW} showed that when the vector field $V$ is the gradient of some smooth function $-f$ on $M$, namely, $V= - \n f$, if the Bakry--Emery Ricci curvature ${\rm Ric}_f$ has a positive lower bound and $|f|$ is bounded, then $M$ is compact. Limoncu \cite{Lim1} obtained that the complete Riemannian manifold $M$ is compact when the $m$-Bakry--Emery Ricci curvature ${\rm Ric}^m_{V}$ has a positive lower bound. Afterwards, Soylu \cite{Soy} and Wang \cite{Wang} verified the Cheeger--Gromov--Taylor type compactness theorem via Bakry--Emery Ricci curvature ${\rm Ric}_f$ and $m$-Bakry--Emery Ricci curvature ${\rm Ric}^m_{f}$ respectively.

In this note, we shall follow the idea of Wan \cite{Wan}, by using the excess function and $V$-Laplacian comparison theorem (see Theorem 3 in \cite{CJQ}), we derive an extension of Bonnet--Myers' type theorem with the Bakry--Emery Ricci curvature.
 
 \begin{thm}\label{thm-2}

Let $(M^n, g)$ be a complete Riemannian manifold, $V$ a smooth vector field on $M$, and $h: [0, +\infty) \to (0, +\infty)$ a continous function. Let $r(x)=d(p, x)$ be the distance function from $p \in M$. Assume $\la V, \n r \ra \leq C_1$ along a minimal geodesic from every point $\widetilde{p} \in M$, here $C_1$ is a constant. Suppose 
\[
{\rm Ric}_V(x) \geq C_2 h(r(x)),
\]
where  $C_2$ is a positive constant depending only on $h, n$ and $C_1$.
Then $M$ is compact. Here $C_2$ can be chosen as $\left( \frac{n-1}{\d} + 2C_1 \right) \cdot \frac{1}{\int_{\d}^{+\infty}h(s)ds}+ \e$ ($ \e, \d $ are arbitrary positive constants).

\end{thm}

\begin{rem}

If $f$ satisfies $\int_{\d}^{+\infty}h(s)ds =+\infty$, then the constant $C_2$ can be chosen as arbitrary positive real number. And by the Gaussian shrinking Ricci soliton, we know that the condition on $\la V, \n r \ra$ can not be weakened to "from a fixed point $\widetilde{p}\in M$".

\end{rem}

\begin{rem}

When $V = - \n f$ for some smooth function $f$ on $M$, if we choose $h\equiv cosntant$, then Theorem \ref{thm-2} becomes Theorem 1.3 in \cite{Wu}.

\end{rem}

Choosing $h(x)=\frac{1}{(r(x)+r_0)^k}$, where $k\in \mathbb{R}$ and $r_0 $ is a positive constant, then Theorem \ref{thm-2} implies 

\begin{cor}\label{cor-2}

Let $(M^n, g)$ be a complete Riemannian manifold and $V$ a smooth vector field on $M$. Suppose that $\la V, \n r \ra \leq C_1$ along a minimal geodesic from every point $\widetilde{p} \in M$, here $C_1$ is a constant. If there exist $p \in M, k\in \mathbb{R}$, and $r_0 >0$, such that 
\[
{\rm Ric}_{V}(x) \geq \frac{C(n, k, r_0, C_1)}{(r(x)+r_0)^k}
\]
for all $r(x)\geq 0$, where $r(x)=d(p, x)$ is the distance function and $C(n, k, r_0, C_1)$ is a constant depending only on $n, k, r_0, C_1$.  Then $M$ is compact. Here when $k>2$, $C(n, k, r_0, C_1)$ can be chosen as $[2C_1r_0+(n-1)(k-2)]\cdot \frac{(k-1)^k}{(k-2)^{k-1}}\cdot r_0^{k-2} + \e$, when $1< k \leq 2$,  $C(n, k, r_0, C_1)$ can be chosen as $\left( 2C_1+\frac{n-1}{\d} \right)(k-1)(r_0+\d)^{k-1} + \e$($ \e, \d $ are arbitrary positive constants), when $k\leq 1$, $C(n, k, r_0, C_1)$ can be chosen as arbitrary positive real number.

\end{cor}

\begin{rem}

When $V \equiv 0$, then $C_1$ can chosen as zero, ${\rm Ric}_V$ is just the usual Ricci curvature, and  Corollary \ref{cor-2} becomes Theorem 0.1 in \cite{Wan}, which is an extension of Calabi and Cheeger--Gromov--Taylor's theorems.

\end{rem}

%\begin{rem}

%When k=2, Tadano (2016) showed that if the Bakry-Emery curvature ${\rm Ric}_V$  is bounded from lower by $\frac{(n-1)(\frac{1}{4}+v^2)}{r^2(x)}$ and the vector field $V$ satisfies $|V|(x)\leq \frac{(n-1)c}{r(x)}$ ($0\leq c < v^2$ is a constant), then $M$ is compact. In this case, the condition on $V$ in Corollary \ref{cor-2} is weaker than the one of his result. 

%\end{rem}

In particular, if we take $h \equiv constant$ in Theorem \ref{thm-2}, then by Remark 1, we obtain

\begin{cor}\label{cor-3}

Let $(M^n, g)$ be a complete Riemannian manifold and $V$ a smooth vector field on $M$. Let $r(x)=d(p, x)$ be the distance function from $p \in M$. If $V$ satisfies $\la V, \n r \ra \leq C_1$ along a minimal geodesic from every point $\widetilde{p} \in M$, here $C_1$ is a constant. Suppose that
\[
{\rm Ric}_{V}(x) \geq (n-1)K,
\]
where $K$ is a positive constant. Then $M$ is compact.

\end{cor}

\begin{rem}

Fern\'andez-L\'opez and Garc\'ia-R\'io \cite{FG} proved that if the norm of $V$ is bounded and the Bakry--Emery Ricci curvature has a positive lower bound, then $M$ is compact (see Theorem 1 in \cite{FG}).  Clearly, the condition on the vector field $V$ in Corollary \ref{cor-3} is weaker than theirs.

\end{rem}

When the vector field $V = - \n f$ for some smooth real valued function $f$ on $M$, applying Theorem 1.1 b) in \cite{WW} instead of 
 the $V$-Laplacian comparison theorem in the proof of Theorem \ref{thm-2}, we can derive the following extension of Bonnet--Myers' type theorem. Since the proof is quite similar to the one of Theorem \ref{thm-2}, we omit the proof here. 
 
  \begin{thm}\label{thm-3}

Let $(M^n, g)$ be a complete Riemannian manifold, $f$ a smooth real valued function on $M$, and $h: [0, +\infty) \to (0, +\infty)$ a continous function. Let $r(x)=d(p, x)$ be the distance function from $p \in M$. Suppose $ |f| \leq C_3$ and 
\[
{\rm Ric}_f(x) \geq C_4 h(r(x)),
\]
where  $C_3\geq 0$ is a constant and $C_4$ is a positive constant depending only on $h, n$ and $C_3$.
Then $M$ is compact. Here $C_4$ can be chosen as $\frac{n+4C_3-1}{\d}   \cdot \frac{1}{\int_{\d}^{+\infty}h(s)ds}+ \e$ ($ \e, \d $ are arbitrary positive constants).

\end{thm}

\begin{rem}

If we take $h \equiv constant$, then $C_4$ can be chosen as arbitrary positive real number. Thus Theorem \ref{thm-3} recovers the Myers'  theorem in \cite{WW}.

\end{rem}

 \begin{cor}\label{cor-4}

Let $(M^n, g)$ be a complete Riemannian manifold and $f$ a smooth real valued function on $M$. Suppose $ |f| \leq C_3$, here  $C_3$ is a nonegative constant. If there exist $p \in M, k\geq 2$, and $r_0 >0$, such that 
\[
{\rm Ric}_{f}(x) \geq \frac{C(n, k, r_0, C_3)}{(r(x)+r_0)^k}
\]
for all $r(x)\geq 0$, where $r(x)=d(p, x) $ is the distance function and $C(n, k, r_0, C_3)$ is a constant depending only on $n, k, r_0, C_3$, and   Then $M$ is compact. Here when $k>2$, $C(n, k, r_0, C_3)$ can be chosen as $(n+4C_3-1)\cdot \frac{(k-1)^k}{(k-2)^{k-2}}\cdot r_0^{k-2} + \e$, when $ k = 2$,  $C(n, k, r_0, C_3)$ can be chosen as $(n+4C_3-1)(1+\frac{r_0}{\d}) + \e$($ \e, \d $ are arbitrary positive constants).

\end{cor}
 
 \begin{rem}
 
 Soylu \cite{Soy} proved a Cheeger--Gromov--Taylor type theorem via the Bakry--Emery Ricci curvature ${\rm Ric}_f$ with bounded $f$, so the case of $k<2$ is covered by his result, we only need to consider the case when $k\geq 2$.
 
 \end{rem}
 
 For  complete Riemannian manifolds with the $m$-Bakery--Emery Ricci curvature ${\rm Ric}^m_{V}$, we can also obtain the Bonnet--Myers' type theorem by using a similar method as in the proof of Theorem \ref{thm-2}.
 
 \begin{thm}\label{thm-1}

Let $(M^n, g)$ be a complete Riemannian manifold, $V$ a smooth vector field on $M$ and $h: [0, +\infty) \to (0, +\infty)$  a continuous function. If there is a positive constant $C_5$ depending only on $h$ and $n$, such that
\[
{\rm Ric}^m_{V}(x) \geq C_5h(r(x)),
\]
where $r(x)= d(p, x)$ is the distance function from $p\in M$. Then $M$ is compact. Here $C_5$ can be chosen as $\frac{n-1}{\d}\cdot \frac{1}{\int_{\d}^{+\infty}h(s)ds} + \e$ ($ \e, \d $ are arbitrary positive constants).

\end{thm}

 \begin{rem}

If $h$ satisfies $\int_{\d}^{+\infty}h(s)ds =+\infty$, then the constant $C_5$ can be chosen as arbitrary positive real number. Hence if we choose $h$ as a constant, then Theorem \ref{thm-1} is just Theorem 1.2 in \cite{Lim1}.

\end{rem}

\begin{cor}\label{cor-1}

Let $(M^n, g)$ be a complete Riemannian manifold and $V$ a smooth vector field on $M$. If there exist $p \in M, k\in \mathbb{R}$, and $r_0 >0$, such that 
\[
{\rm Ric}^m_{V}(x) \geq \frac{C(n, k, r_0)}{(r(x)+r_0)^k}
\]
for all $r(x)\geq 0$, where $r(x)= d(p, x)$ is the distance function and $C(n, k, r_0)$ is a constant depending only on $n, k, r_0$, then $M$ is compact. Here when $k>2$, $C(n, k, r_0)$ can be chosen as $(n-1)\cdot \frac{(k-1)^k}{(k-2)^{k-2}}\cdot r_0^{k-2} + \e$, when $1< k \leq 2$,  $C(n, k, r_0)$ can be chosen as 
 $\frac{n-1}{\d}\cdot (k-1)\cdot(r_0+\d)^{k-1} + \e$($ \e, \d $ are arbitrary positive constants),
when $k\leq 1$, $C(n, k, r_0)$ can be chosen as arbitrary positive real number.

\end{cor}

\begin{rem}

When $V= - \n f$ for some smooth function $f$ on $M$, Wang \cite{Wang}  proved a Cheeger--Gromov--Taylor type compactness theorem via $m$-Bakry--Emery curvature ${\rm Ric}^m_{f}$. 

\end{rem}

\vskip24pt

\section{Proof of Theorem \ref{thm-2} and Theorem \ref{thm-1}}

Denote the $V$-Laplacian operator $\D_V:=\D+\la V, \n\cdot \ra$. 

\vskip 12pt

\noindent{\it Proof of Theorem \ref{thm-2}}  \quad
Suppose that $M$ is  complete noncompact. Then for any $p\in M$, there is a ray $\s(t)$, such that $\s(0) = p$.

Since $r(x)=d(p, x)$ is smooth outside the cut-locus of $p$, applying the Bochner formula (2.2) in \cite{CQ}, we have
\begin{equation}\label{eqn-Boch}
|{\rm Hess}(r)|^2 + {\rm Ric}_V (\n r, \n r) + \la  \n \D_V r, \n r \ra = \frac{1}{2} \D_V |\n r|^2 = 0.
\end{equation}

Set $\varphi_V(t):= (\D_V r) \cdot \s(t)$. Computing both sides of (\ref{eqn-Boch}) along $\s(t)$, 
\begin{equation*}
|{\rm Hess}(r(\s(t))|^2 = - {\rm Ric}_V(\s'(t), \s'(t)) - \varphi'_V(t).
\end{equation*}
$\forall \d>0$, integrating the above equality over the interval $[\d, t]$,
\begin{equation}\label{eqn-Boch5}
\int_{\d}^{t} |{\rm Hess}(r(\s(s))|^2 ds \leq \varphi_V(\d) - \varphi_V(t) - \int_{\d}^{t} {\rm Ric}_{V} (\s'(s), \s'(s) )ds.
\end{equation}
Since ${\rm Ric}_V(x) \geq C_2 h(r(x))>0$ and $\la V, \n r \ra \leq C_1$, by the V-Laplacian comparison theorem  (see Theorem 3 in \cite{CJQ}), 
\[
\varphi_V(t)=(\D_V r)\circ \s(t) \leq \frac{n-1}{t}+C_1.
\]
Consider the excess function
\[
e(x) := d(p, x) + d(\s(i), x) - i. 
\]
By the triangle inequality, we have
\[
e(x)\geq 0 \quad {\rm and} \quad e(\s(t)) \equiv 0 \quad {\rm for } \quad 0\leq t\leq i.
\]
Hence 
\[
\D_V d(p, \s(t)) + \D_V d(\s(i), \s(t)) = \D_V e(\s(t)) \geq 0.
\]
It follows that 
\[
\varphi_V(t) = (\D_V r) \circ \s(t) = \D_V d(\s(0), \s(t)) \geq -\D_V d(\s(i), \s(t)) \geq -\frac{n-1}{i-t}-C_1.
\]
%Then using the excess function as in the proof of Theorem \ref{thm-1}, we can conclude that 
%\[
%\varphi_V(t) \geq -\frac{n-1}{i-t} - a.
%\]
Letting $i \to +\infty$, we obtain $\varphi_V(t) \geq -C_1$. Therefore we arrive at 
\begin{equation}\label{eqn-Lap2}
-C_1\leq \varphi_V(t) \leq \frac{n-1}{t}+C_1.
\end{equation}
From (\ref{eqn-Boch5}) and (\ref{eqn-Lap2}), we derive
\[
0\leq \int_{\d}^{t} |{\rm Hess}(r(\s(s))|^2 ds \leq \frac{n-1}{\d} + 2C_1 - \int_{\d}^t C_2 h(s)ds.
\]
Letting $t \to +\infty$, we get
\[
C_2\int_{\d}^{+\infty} h(s) ds \leq \frac{n-1}{\d}+2C_1.
\]
Choosing $C_2=\left( \frac{n-1}{\d} + 2C_1 \right) \cdot \frac{1}{\int_{\d}^{+\infty}h(s)ds}+ \e$ ($\forall \e > 0$ is a constant), this yields the contradiction. Hence $M$ must be compact.
\qed

\vskip12pt

\noindent{\it Proof of Theorem \ref{thm-1}}  \quad
Suppose that $M$ is complete noncompact. Then for any $p\in M$, there is a ray $\s(t)$, such that $\s(0) = p$. 

By the elementary inequality 
\[
(a+b)^2 \geq \frac{a^2}{1+\a} - \frac{b^2}{\a}, \quad \forall \a,
\]
we get 
\begin{equation}\label{eqn-Hess}\aligned
|{\rm Hess}(r)|^2 \geq & \frac{1}{n} (\D r)^2 = \frac{1}{n} (\D_V r - \la V, \n r \ra)^2 \\
 \geq & \frac{1}{n} \left( \frac{1}{\frac{m}{n}} (\D_V r)^2 - \frac{1}{\frac{m}{n}-1} \la V, \n r \ra^2  \right) 
 = \frac{1}{m} (\D_V r)^2 - \frac{1}{m-n} \la V, \n r \ra^2.
\endaligned
\end{equation}
Substituting (\ref{eqn-Hess} ) into (\ref{eqn-Boch}), 
\begin{equation}\label{eqn-Boch2}
\frac{1}{m} (\D_V r)^2 + {\rm Ric}^m_{V} (\n r, \n r) + \frac{\p}{\p r} (\D_V r) \leq 0.
\end{equation}
Set $\varphi_V (s) := (\D_V r) \circ \s(t)$. Computing both sides of (\ref{eqn-Boch2}) along $\s(t)$ gives
\begin{equation}\label{eqn-Boch3}
\frac{1}{m} (\varphi_V (t))^2 + {\rm Ric}^m_{V} (\s'(t), \s'(t) ) + \varphi'_V (t) \leq 0.
\end{equation}
For any $\d>0$, integrating (\ref{eqn-Boch3}) over the interval $[\d, t]$, we obtain
\begin{equation}\label{eqn-Boch4}
\frac{1}{m}\int_{\d}^{t} (\varphi_V (s))^2 ds \leq \varphi_V(\d) - \varphi_V(t) - \int_{\d}^{t} {\rm Ric}^m_{V} (\s'(s), \s'(s) )ds.
\end{equation}
Since ${\rm Ric}^m_{V}(x) \geq C_5h(r(x)) > 0$,  by Theorem 4.2 in \cite{BQ}, $\varphi_V(t)=(\D_V r)\circ \s(t) \leq \frac{n-1}{t}$.
Then using the excess function as in the proof of Theorem \ref{thm-2}, we can conclude that 
\[
\varphi_V(t) \geq -\frac{n-1}{i-t}.
\]
Letting $i \to \infty$, we get $\varphi_V(t) \geq 0$. Therefore we obtain 
\begin{equation}\label{eqn-Lap}
0\leq \varphi_V(t) \leq \frac{n-1}{t}.
\end{equation}
Combing (\ref{eqn-Boch4}) with (\ref{eqn-Lap}), it follows that
\[
0\leq \frac{1}{m}\int_{\d}^t (\varphi_V(s))^2 ds \leq \frac{n-1}{\d} - \int_{\d}^t C_5h(s) ds. 
\] 
Letting $t \to +\infty$, we get
\[
C_5\int_{\d}^{+\infty} h(s) ds \leq \frac{n-1}{\d}. 
\]
Choosing $C_5 = \frac{n-1}{\d}\cdot \frac{1}{\int_{\d}^{+\infty}h(s)ds}+\e$ ($\forall \e > 0$ is a constant), it is a contradiction. Hence $M$ must be compact.
\qed

\vskip12pt

\noindent{\bf Acknowledgement}  This work is partially supported by NSFC (Nos. 11771339, 11571259), Fundamental Research Funds for the Central Universities (No. 2042019kf0198) and the Youth Talent Training Program of Wuhan University. The author thanks Dr. Jianming Wan  for providing him with the work of \cite{Wan} and for  helpful discussions.  He  would like to express his gratitude to Professor Tobias H. Colding for his invitation, to MIT for their hospitality.

\end{document}